\documentclass[a4paper]{amsart}
\usepackage{amsfonts}
\theoremstyle{plain}
\begingroup
\newtheorem{theorem}{Theorem}[section]
\newtheorem{lemma}[theorem]{Lemma}

\newtheorem{remark}[theorem]{Remark}
\newtheorem{definition}[theorem]{Definition}
\theoremstyle{definition}
\theoremstyle{remark}
\numberwithin{equation}{section}
\newcommand{\e}{\varepsilon}
\newcommand{\Om}{\Omega}

\newcommand{\ci}{\stackrel{i}{\to}}
\newcommand{\m}{{\bf m}}
\newcommand{\dx}{\,dx}
\newcommand{\ho}{{\rm hom}}

\newcommand{\wup}{{\rm W}^{1,p}}
\newcommand{\orn}{(\Omega;\R^N)}

\newcommand{\msim}{{\rm M}^{N\times N}_{\rm sym}} 
\newcommand{\lm}{{\mathcal L}^{N}}
\newcommand{\hn}{{\mathcal H}^{N-1}} 

\newcommand{\F}{{\mathcal F}} 
\newcommand{\E}{{\mathcal E}}
\newcommand{\R}{{\mathbb R}}
\newcommand{\Op}{{\mathcal O}} 
\newcommand{\tr}{{\rm tr}\,}

\newcommand{\salt}{\noalign{\vskip .2truecm}}
\newcommand{\parent}[3]{\left #1 {#3} \right #2} 
\newcommand{\graffe}[1]{\parent \{ \}{#1}} 
\newcommand{\res}{\mathop{\hbox{\vrule height 7pt width .5pt depth 0pt
\vrule height .5pt width 6pt depth 0pt}}\nolimits}

\newcommand{\sbd}{{\rm SBD}(\Om)}
\newcommand{\bd}{{\rm BD}(\Om)}

\newcommand{\bom}{{\mathcal B}(\Om)} 
\newcommand{\opom}{{\mathcal O}(\Om)}
\newcommand{\omreg}{{\mathcal O}_{\infty}(\Om)} 
\newcommand{\intcauchy}{\mskip 5mu -\mskip -18mu \int}

\newcommand{\norm}[1]{\left\Vert {#1} \right\Vert} 
\title{A note on the integral representation of functionals in the space 
$\sbd$}
\author[Fran\c{c}ois Ebobisse]{Fran\c {c}ois Ebobisse}
\address[Fran\c{c}ois Ebobisse]{S.I.S.S.A., Via Beirut 2-4, 34014, Trieste, 
Italy}
\email[F. Ebobisse]{ebobisse@sissa.it}
\author[Rodica Toader]{Rodica Toader}
\address[Rodica Toader]{S.I.S.S.A., Via Beirut 2-4, 34014, Trieste, Italy}
\email[R. Toader]{toader@sissa.it}

\begin{document}
\vskip .2truecm
\begin{abstract}
\small{In this paper we study the integral representation in the space ${\rm SBD}$ of 
special functions with bounded deformation of some $L^1$-norm lower 
semicontinuous functionals invariant with respect to rigid motions.
 }
\end{abstract}
\maketitle

\section{Introduction}

Several phenomena in phase transition, fracture mechanics, liquid crystals, 
can be mo\-delled as energy minimization problems where the natural energy 
has both volume and surface terms. In many cases the energy functional is 
obtained as a limit of approximating functionals and some of its properties 
can be deduced from the approximation process.

A basic step is then to obtain, starting from these properties, an integral 
representation of the energy. We consider here this problem for local 
functionals 
$\F$ defined on the space ${\rm BD}$ of functions with bounded deformation, 
which are lower semicontinuous with respect to the $L^1$-topology, satisfy 
linear growth and coercivity conditions, as set functions are 
(restrictions of) Radon measures, and are invariant with respect to rigid 
motions.
In order to identify the volume and the surface densities we follow the 
global method for relaxation introduced by Bouchitt\'e, Fonseca and 
Mascarenhas in \cite{BFM} for functionals defined on 
the space ${\rm BV}$ of functions with bounded variation, 
which is characterized by the identification of 
both bulk and surface densities from a local Dirichlet problem. This kind of 
approach has already been used in some other contexts, as, for instance, 
homogenization, where the homogenized density is obtained from a Dirichlet 
problem in the cell.

An example of functional in the class we consider is given by 
 the relaxed functional $\overline F$ of the bulk energy
\begin{equation}\label{bft}
F(u):=\left\{\begin{array}{ll}\displaystyle\int_\Om f(Eu(x))\,dx & \hbox{if }u\in 
W^{1,1}(\Om,\R^N)\\
+\infty & \hbox{otherwise}
\end{array}
\right.\,
\end{equation}
with respect to the $L^1$-norm topology, where $\Om$ is a bounded open subset 
of $\R^N$ and $f$ is a Borel function satisfying standard linear growth 
assumptions. The integral representation of $\overline F$ on $\bd$ was studied by 
Barroso, Fonseca and Toader in \cite{BFT}, where the global method was applied in order to derive
the surface density, while the volume density was obtained by a direct proof using the explicit form of 
the functional $F$. 

In this paper, the bulk density is deduced from the global method and  
the {\it approximate differentiability} of ${\rm BD}$ functions proved by 
Ambrosio, Coscia and Dal Maso in \cite{ACDM}, while the surface density is 
obtained exactly as in \cite{BFT}. 
Note that both our result and the one in \cite {BFT} are valid for functions 
in  $\sbd$, i.e. integrable functions $u$ for which the Cantor part $E^cu$ of 
the measure $Eu$ vanishes. An integral representation in all the space 
$\bd$ would 
require more information on $E^cu$, since  
the only property that $E^cu$ vanishes on $\hn$-$\sigma$ finite Borel 
subsets, proved in \cite{ACDM}, is not sufficient. We recall in Section 2 some 
useful properties of ${\rm BD}$ functions.

In Section 3 we prove the integral representation theorem 
(Theorem \ref{integ}) and give an example showing why we assume the invariance 
with respect to rigid motions. In the last section we apply 
Theorem~\ref{integ} to obtain the integral representation in $\sbd$ of 
some $\Gamma$-limits arising in the 
homogenization of multi-dimensional structures recently studied in the context 
of linear elasticity and perfect plasticity by Ansini and Ebobisse in 
\cite{AE}, following the measure-theoretic approach introduced by Ansini, 
Braides and Chiad\`o Piat in \cite{ABC}. 

\section{Notation and preliminaries}\label{prel}

Let $N\geq 1$ be an integer. 
We denote by $M^{N\times N}$ the space of $N\times N$ matrices and by 
$\msim $ the subspace of symmetric matrices in $M^{N\times N}$. For any 
$\xi\in M^{N\times N}$, $\xi ^T$ is the transposition of $\xi$. 
Given $u,\,v\in\R ^N$, $u\otimes v$ and $u\odot v:=(u\otimes v+v\otimes u)/2$ 
denote the tensor and symmetric products of $u$ and $v$, respectively. 
We use the standard notation, $\lm $ and $\hn $ to denote respectively the 
Lebesgue and $(N-1)$-dimensional Hausdorff measures. 

Let $\Om $ be a bounded open subset of $\R ^N$; we denote by $\bom$, $\opom$ 
and $\omreg$ the family of Borel, open and open subsets of $\Om$ with 
Lipschitz boundary, respectively.
For any $x\in \Om $ and $\rho>0$, we denote by $B(x,\rho )$ the open ball of 
$\R^N $ centered at
$x$ with radius $\rho $, by $Q(x,\rho)$ the cube of centre $x$ and sidelength 
$\rho $, while $Q_\nu (x,\rho )$ is the cube with two its faces perpendicular 
to the unit vector $\nu$. When $x=0$ and $\rho=1$ we simply write $B$ and $Q$. 
If $\mu $ is a Radon measure, we denote by $|\mu |$ its total variation.

\begin{definition}\label{bd}
A function $u:\Om \to \R^N $ is with {\em bounded deformation} in $\Om $ 
if $u\in L^1(\Om ,\R^N)$ and 
$Eu:=(Du+Du^T)/2\in M_b\bigl(\Om ,\msim\bigr)$, 
where $Du$ is the distributional gradient of $u$ and 
$M_b\bigl(\Om,\msim\bigl)$ 
is the space of $\msim$-valued
Radon measures with finite total variation in $\Om $.
\end{definition}
The space $\bd$ of functions with bounded deformation in $\Om$, introduced in 
\cite{MSC}, has been widely
studied, for instance by Anzellotti-Giaquinta \cite{AG}, Kohn \cite{KH},
Suquet \cite{SUQ}, and Temam \cite{TEM}.
It is a Banach space when equipped with the norm
$$
\norm{u}_{BD(\Om )}:=\norm{u}_{L^1(\Om ,\R ^N)}+\vert Eu\vert (\Om).
$$

It is sometimes convenient to consider also the distance between two functions 
$u$, $v\in\bd$ given by
$$
d(u,v):= \norm{u-v}_{L^1(\Om ,\R ^N)}+|\,|Eu|(\Om)-|Eu|(\Om)|.
$$
The topology induced by this distance in $\bd$ is called 
{\it intermediate topology}. We denote by $\ci$ the convergence with respect 
to this topology. 
It is well known (see Temam \cite{TEM}) that the trace operator 
$\tr:\bd\to L^1(\Gamma ,\R^N)$ is continuous when $\bd$ is equipped with the 
intermediate topology.

Whenever the open set $\Om $ is assumed to be connected, the kernel of the 
operator $E$ is the class of {\it rigid motions} denoted here by 
${\mathcal R}$, and composed of affine maps of the form $Mx+b$, where $M$ is a 
skew-symmetric $N\times N$ matrix and $b\in\R^N$. Therefore ${\mathcal R}$ is 
closed and finite-dimensional.

Fine properties of ${\rm BD}$ functions were studied, for instance, in 
\cite{ACDM}, \cite{BCDM} and \cite{KH}. We recall that if $u\in\bd$, then the 
jump set $J_u$ of $u$ is a countably $(\hn,n-1)$-rectifiable Borel set and the 
following decomposition of the measure $Eu$ holds
\begin{equation}\label{decomp}
Eu=\E u\lm + E^su=\E u\lm +([u]\odot\nu_u)\hn\res J_u +E^cu\,,
\end{equation}
where $[u]:= u^+-u^-$, $u^+$ and $u^-$ are the {\it one-sided Lebesgue limits} of $u$ with 
respect to the measure theoretic normal $\nu_u$ of $J_u$, 
$\E u$ is the density of the absolutely continuous part of $Eu$ with respect 
to $\lm$, $E^su$ is the singular part, and $E^cu$ is the {\it Cantor part} 
and vanishes on the Borel sets that are $\sigma$-finite with respect to $\hn$ 
(see \cite{ACDM}). 

Moreover, the following theorem on the approximate differentiability of 
${\rm BD}$ functions was proved in \cite{ACDM}.

\begin{theorem}\label{approx}
Let $\Om $ be a bounded open set in $\R ^N$ with Lipschitz boundary. Let 
$u\in\bd$. Then for $\lm$ almost every $x\in\Om $ there exists an $N\times N$ matrix 
$\nabla u(x)$ such that
\begin{equation}\label{apdif}
\lim_{\rho\to0}\frac{1}{\rho ^N}\int_{Q(x,\rho)}\frac{|u(y)-u(x)-
\nabla u(x)(y-x)|}{\rho}dy=0\,,
\end{equation}
and 
\begin{equation}\label{symapp}
\lim_{\rho\to0}\frac{1}{\rho ^N}\int_{Q(x,\rho)}
\frac{|(u(y)-u(x)-\E u(x)(y-x),y-x)|}{|y-x|^2}dy=0
\end{equation}
for $\lm$ almost every $x\in\Om $.
\end{theorem} 
In particular, by (\ref{apdif}) $u$ is approximately differentiable 
$\lm$ almost everywhere in 
$\Om$ and the function $\nabla u$ satisfies the weak $L^1$ estimate
$$
\lm (\{x\in\Om\mbox{ : } |\nabla u(x)|>t\}\leq
\frac{C(N,\Om)}{t}\norm{u}_{\bd}\quad\forall t>0,
$$
where $C(N,\Om)$ is a positive constant depending only on $N$ and $\Om $.

From (\ref{symapp}) and (\ref{apdif}) one can easily see that 
\begin{equation}\label{symappdiff}
\E u(x)=(\nabla u(x)+\nabla u(x)^T)/2\quad\mbox{for }\lm \mbox{-a.e. }x\in\Om.
\end{equation}

\noindent Analogously to the space ${\rm SBV}$ introduced by De Giorgi and 
Ambrosio in \cite{DA}, the space ${\rm SBD}$ was defined in \cite{BCDM}.
\begin{definition}\label{sbd} 
The space $\sbd$ of {\em special functions with bounded deformation}, is the 
space of functions $u\in\bd$ such that the measure $E^cu$ in (\ref{decomp}) 
is zero.
\end{definition}

\section{Main result}\label{main}

Let $\F:\bd\times\omreg\to[0,+\infty]$ be a functional satisfying the 
properties mentioned in the introduction, more precisely,
\begin{enumerate}
\item[{\rm (}1{\rm )}] $\F(\cdot,A)$ is $L^1(A,\R^N)$ lower semicontinuous;
\item[{\rm (}2{\rm )}] for every $u\in\bd$, 
\begin{equation}\label{2.3.9}
\frac{1}{C}|Eu|(A)\leq \F(u,A)\leq C(\lm(A)+|Eu|(A));
\end{equation}
\item[{\rm (}3{\rm )}] $\F(u,\cdot)$ is the restriction to $\omreg$ of a Radon 
measure;
\item[{\rm (}4{\rm )}] $\F(u+R)=\F(u)$ for every $u\in\bd$ and every 
rigid motion $R$.
\end{enumerate}
Since the properties (2) and (3) give the absolute continuity of 
${\F}(u, \cdot)$ with respect to the measure $\mu:=\lm +|E^su|$, in order to 
obtain the integral representation of $\F$, we need only 
to identify the volume and the surface densities whenever $u\in\sbd$.
To do this we define, as in \cite{BFM}, see also \cite{BFT}, 
for every $u\in\bd$ and every $A\in \omreg$
$$
\m(u,A):=\inf\{{\F}(v, A) : v\in\bd\,,\quad v|_{\partial A}= u|_{\partial A}\}.
$$
The basic idea of the global method in \cite{BFM} consists in comparing the 
asymptotic behaviours of $\m(u, Q(x_0,\e))$ and ${\F}(u,Q(x_0,\e) )$ with 
respect to $\mu(Q(x_0,\e))$ as $\e\to 0^+$, and to show via a blow-up argument 
that, the volume and surface densities are obtained from a local Dirichlet 
problem (see Lemma~\ref{3.5bfm}).\\
We shall use the following lemmas, similar to Lemmas~3.1 and~3.5 
in~\cite{BFM} for ${\rm BV}$-functions, proved in the case of 
${\rm BD}$-functions in \cite[ Lemmas~3.10,~3.12]{BFT}.
\begin{lemma}\label{m}{\rm (\cite[Lemma 3.10]{BFT})} 
There exists a positive constant $C$ such that for any $u_1$, $u_2\in \bd$ and 
any $A\in{\Op}_{\infty}(\Om)$ we have
$$ 
|\m(u_1,A)-\m(u_2,A)|\leq C\int_{\partial A}|\tr(u_1-u_2)(x)|d\hn(x).
$$
\end{lemma}
\begin{lemma}\label{3.5bfm}{\rm (\cite[Lemma 3.12]{BFT})}
If $\F$ satisfies conditions (1)-(3)  then
$$ 
\lim_{\e\to0}\frac{\F(u,Q_{\nu}(x_0,\e))}{\mu(Q_{\nu}(x_0,\e))}=
\lim_{\e\to0}\frac{\m(u,Q_{\nu}(x_0,\e))}{\mu(Q_{\nu}(x_0,\e))}\;\;
\mu\;a.e.\;x_0\in\Omega\hbox{ and for all }\nu\in S^{N-1}.
$$
\end{lemma}
%
We prove now the integral representation result.
\begin{theorem}\label{integ}
Let $\F:\bd\times\omreg\to[0,+\infty]$ be a functional satisfying 
properties (1)-(4). Then for every $u\in\sbd$ and $A\in\omreg$ we have
\begin{equation}\label{int1}
\F(u,A)=\int_A f(x,\E u)dx +\int_{J(u)\cap A} g(x,[u],\nu)d\hn\,,
\end{equation}
where
\begin{eqnarray}
f(x_0,\xi)&:=&\limsup_{\e\to0}\frac{\m(\xi(\cdot-x_0),Q(x_0,\e))}{\e^N}
\label{f}
\\
g(x_0,\lambda,\nu)&:=&\limsup_{\e\to0}\frac{\m(u_{\lambda,\nu}
(\cdot-x_0),Q_\nu(x_0,\e))}{\e^{N-1}}\label{g}
\end{eqnarray}
for all $x_0\in\Om$, $\lambda\in\R^N$, $\xi\in \msim$, 
$\nu\in S^{N-1}$, and where 
$$
u_{\lambda,\nu}(y):=\left\{
\begin{array}{ll}\lambda & \hbox{if } y\cdot\nu>0\\
0 & \hbox{otherwise.}
\end{array}\right.
$$
\end{theorem}

We use the same notation for $\F(u,\cdot)$ and its extension to the Borel 
subsets of $\Om$.

\begin{proof} (i) {\it The volume part}.  
Let $u\in \sbd$ and choose $x_0\in \Om$ such that
\begin{eqnarray}
&& \frac{d\F(u,\cdot)}{d\lm}(x_0)=\lim_{\e\to0}\frac{\F(u,Q(x_0,\e))}{\e^N}\;
\hbox{ exists and is finite,}\label{4.1}\\ 
&& \lim_{\e\to0}\frac{1}{\e^{N+1}}
\int_{Q(x_0,\e)}|u(x)-u(x_0)-\nabla u(x_0)(x-x_0)|dx=0,\label{4.2}\\ 
&& \lim_{\e\to0}\frac{1}{\e^N}|Eu|(Q(x_0,\e))=|\E u(x_0)|,\label{4.3}\\
&& \lim_{\e\to0}\frac{1}{\e^N}|E^su|(Q(x_0,\e))=0,\label{4.4}\\
&& \frac{d\F(u,\cdot)}{d\lm}(x_0)=\lim_{\e\to0}\frac{\m(u,Q(x_0,\e))}{\e^N}.
\label{dfdlm}
\end{eqnarray}
Let, for every $y\in Q$,
$$ 
u_{\e}(y):=\frac{u(x_0+\e y)-u(x_0)}{\e}\qquad \hbox{and}\qquad
u_0(y):=\nabla u(x_0)y.
$$ 
By (\ref{4.2}) the functions $u_\e$ converge to $u_0$ in $L^1(Q,\R^N)$. 
Moreover, 
$$
|Eu_\e|(Q)\to|Eu_0|(Q).
$$
Indeed, by definition
\begin{eqnarray}
\nonumber |Eu_\e|(Q)&=&\sup_{\stackrel{\phi\in C_0^1\bigl(Q,\,\msim\bigr)} 
{\|\phi\|_{\infty}\leq 1}}\int_Q\frac{u(x_0+\e y)-u(x_0)}{\e} {\rm div}\,\phi(y)dy\\
\nonumber &=&\sup_{\stackrel{\varphi\in C_0^1\bigl(Q(x_0,\e),\,\msim\bigr)}{
 \|\varphi\|_{\infty}\leq 1}}\frac{1}{\e^N}\int_{Q(x_0,\e)}(u(x)-u(x_0)) 
{\rm div}\,\varphi(x)dx\\
\nonumber &=&\frac{1}{\e^N}|Eu|(Q(x_0,\e)),
\end{eqnarray}
where $\varphi(x):=\phi(\frac{x-x_0}{\e})$. 
Then from (\ref{4.3}) we get $|Eu_\e|(Q)\to|\E u(x_0)|=|Eu_0|(Q)$, where we used 
also the formula (\ref{symappdiff}). This shows that $u_\e\ci u_0$ in ${\rm BD}(Q)$.\\
On the other hand from the continuity of the trace with respect to the 
intermediate topology it follows that 
\begin{eqnarray*} 
&&\int_{\partial Q}|\tr (u_\e(y)-\nabla u(x_0)(y))|d\,\hn(y)\\
&&=\frac{1}{\e^N}
\int_{\partial Q(x_0,\e)}|\tr (u(x)-u(x_0)-\nabla u(x_0)(x-x_0))|d\,\hn(x)\,\to 0.
\end{eqnarray*}
Then by (\ref{dfdlm}), Lemmas \ref{m} and~\ref{3.5bfm} we obtain
\begin{eqnarray}
\nonumber \frac{d\F(u,\cdot)}{d\lm}(x_0) & = &\lim_{\e\to 0}
\frac{\m(u,Q(x_0,\e))}{\e^N}\\
\nonumber & = &  \lim_{\e\to 0}\frac{\m(u(x_0)+\nabla u(x_0)(\cdot-x_0),Q(x_0,\e))}{\e^N}.
\end{eqnarray} 
Now condition (4) with $R(x):=u(x_0)+\frac{\nabla u(x_0)-\nabla u(x_0)^T}{2}(x-x_0)$ implies that
\begin{eqnarray*}
\nonumber \frac{d\F(u,\cdot)}{d\lm}(x_0) & = &\lim_{\e\to 0}\frac{\m(\E u(x_0)(\cdot-x_0),Q(x_0,\e))}{\e^N}\\
& = & f(x_0, \E u(x_0)).
\end{eqnarray*}

(ii) {\it The surface part}.
As in \cite[Proposition~5.1]{BFT}, it can be proved that 
$$\F(u, A\cap J_u)=\int_{J(u)\cap A} g(x,u^+, u^-,\nu)d\hn\,,
$$
where 
$$ 
g(x_0,\lambda,\theta,\nu):=\limsup_{\e\to0}\frac{\m(u_{\lambda,\theta,\nu}
(\cdot-x_0),Q_\nu(x_0,\e))}{\e^{N-1}}
$$
for all $x_0\in\Om$, $\lambda , \theta\in\R^N$, $\nu\in S^{N-1}$, and where 
$$
u_{\lambda,\theta,\nu}(y):=\left\{
\begin{array}{ll}\lambda & \hbox{if } y\cdot\nu>0\\
\theta & \hbox{otherwise.}
\end{array}\right.
$$
Using again condition (4) we obtain
$$
\m(u_{\lambda,\theta,\nu}(\cdot-x_0),Q_\nu(x_0,\e))=\m(u_{\lambda-\theta,\nu}
(\cdot-x_0),Q_\nu(x_0,\e)),
$$
hence (\ref{g}), concluding thus the proof.
\end{proof}
\begin{remark}{\rm
As a particular case the result in \cite{BFT} is recovered, i.e. if $\overline F$ is 
the localized lower semicontinuous envelope of the functional $F$ given by (\ref{bft}), then 
$$
\overline F(u,A)=\int_A SQf({\mathcal E}u(x))\,dx +
\int_{A\cap J_u}
\bigl(SQf\bigr)^\infty \bigl([u]\odot\nu_u(x)\bigr)\,d{\mathcal H}^{N-1}(x)
$$
for every $u\in\sbd$ and every $A\in {\mathcal O}_\infty(\Om )$, where $SQf$ is the 
symmetric quasiconvex envelope of $f$ introduced by Ebobisse in \cite{EBO2}, and 
characterized by 
$$   
SQf(\xi)=\inf\graffe{\intcauchy _A f(\xi+{\mathcal E}\varphi (x))dx\mbox{; }\varphi 
\in W^{1,\infty }_0(A,\R ^N )},
$$
for every $\xi\in\msim $ and for every bounded open subset $A$ of $\R ^N $, and 
$f^\infty$ is the recession function of $f$.
}
\end{remark}

\begin{remark}\label{rigid}
{\rm Note that hypothesis (4) is not a consequence of hypotheses (1)-(3). In fact, without condition (4) of Theorem \ref{integ} we would obtain that
$$
{\F}(u, A)=\int_Af\bigl(x,u(x),\nabla u(x)\bigr)\dx+
\int_{J_u\cap A}g\bigl(x,u^+(x),u^-(x),\nu_u(x)\bigr)\,d{\mathcal H}^{N-1}(x)
$$
for every $ (u,A)\in \sbd\times {\mathcal O}_\infty(\Om )$. In particular, for every 
$u\in W^{1,1}(\Om, \R ^N)$,
$$
{\F}(u, \Om)=\int_\Om f\bigl(x,u(x),\nabla u(x)\bigr)\dx,
$$
which, under some continuity assumption on $\F$ with respect to $u$, 
for instance, assuming that there exists a modulus of continuity $\psi(t)$ satisfying
$$
|\F(u(\cdot-z)+w, z+A)-\F(u,A)|\leq\psi(|w|+|z|)(\lm(A)+|Eu|(A))\,,
$$
for all $(u,A,w,z)\in \bd\times {\mathcal O}_\infty(\Om )\times\R^N\times\R^N$, 
such that $z+A\subset\Om$,
implies that 
 for $\lm$-almost every $x_0\in\Om $ and every $p\in\R ^N$, the function 
$f(x_0,p,\cdot)$ is quasiconvex. 
By (2), 
$$
\frac{1}{C}|\xi+\xi^T|\leq f\bigl(x_0,\xi x_0,\xi\bigr)\leq C\bigl(1+|\xi+\xi^{T}|\bigr).
$$
The following example shows that there exists a rank-one convex function 
$\phi:M^{2\times 2}\to[0,+\infty[$
 which satisfies 
\begin{equation}\label{phi}
\frac{1}{C}|\xi+\xi^T|\leq \phi\bigl(\xi\bigr)\leq C\bigl(1+|\xi+\xi^T|\bigr)\qquad\forall 
\xi\in M^{2\times 2}\,,
\end{equation}
 and which depends also on the  
antisymmetric part of the matrix $\xi$.
Let 
$\xi:=\left(\begin{array}{cc}
a & b\\
c & d
\end{array}\right)$. 
It is enough to define such a function on the matrix 
$\left(\begin{array}{cc}
0 & b\\
c & 0
\end{array}\right)$ and then to add the quantity $|a|+|d|$. Since 
$\left(\begin{array}{cc}
0 & b\\
c & 0
\end{array}\right)=\displaystyle\frac{c+b}{2}\left(\begin{array}{cc}
0 & 1\\
1 & 0
\end{array}\right)+\displaystyle\frac{c-b}{2}\left(\begin{array}{cc}
0 & -1\\
1 & 0
\end{array}\right),$
we look for a function $h(b,c)$ which is separately convex, satisfies 
the linear growth condition  
$$
0\leq h(b,c)\leq C(1+|b+c|)
$$
and which depends also on $c-b$. An example of such a function is the following:
$$
h(x,y):=\left\{\begin{array}{cc}
-(x+y) & \mbox{ if } x+y\leq 0\\
\salt
0 & \mbox{ if }0<x+y\leq 1\mbox{ and }xy\leq 0\\
\salt
xy & \mbox{ if }x,y\geq 0\mbox{ and }\max(x,y)\leq 1\\
\salt
x+y-1 & \mbox{ if }x+y>1\mbox{ and }\max(x,y)> 1.
\end{array}\right.
$$
Therefore, the function $\phi: M^{2\times 2}\to [0,+\infty[$ given by 
$$
\phi(\xi):=h(b,c)+|a|+|d|+|b+c|\qquad\hbox{where }\xi=\left(\begin{array}{cc}
a & b\\
c & d
\end{array}\right)\,,
$$
is rank-one convex, has the linear growth (\ref{phi}), and does 
not depend only on the symmetric part of the matrix $\xi$. 
}
\end{remark}


\section{Application: homogenization of periodic multi-dimensional 
structures}\label{appli}
In  \cite{AE}, the authors studied the asymptotic behaviour of 
functionals of the form 
$$ 
F_\e(u,\Om):=\int_\Om \varphi\Bigl(\frac{x}{\e},\frac{dEu}{d\mu_\e}\Bigr)d\mu_\e 
$$
defined on a particular class of functions
with bounded deformation, denoted by $LD_{\mu _\e }^p(\Om)$, and given by 
the functions $u\in L^p(\Om, \R ^N)$ whose deformation
tensor $Eu$ is an absolutely continuous measure with respect to $\mu_\e $,  
with $p$-summable density $dEu/d\mu_\e $,
where $\mu _\e $ is defined by $\mu _\e (B):=\e ^N\mu (\frac{1}{\e}B)$,
with $\mu $ a fixed $1$-periodic Radon measure and $\varphi$ is a
Borel function $1$-periodic in the first variable. 
Assuming the standard $p$-growth 
condition on $\varphi$ and that the measure $\mu$ is '$p$-homogenizable' 
(see \cite[\S 4]{AE}), the authors 
proved a homogenization theorem (Theorem 5.1). 
Precisely, they proved the
existence of the $\Gamma $-limit $F_{\ho}$ of the functionals $F_\e $ 
with respect to $L^p$-convergence in the
Sobolev space $W^{1,p}(\Om ,\R ^N)$, and with respect to 
$L^1$-convergence in $\bd$. They showed that the
$\Gamma $-limit admits the integral representation
\begin{equation}\label{intw}
F_{\ho}(u,\Om)=\int_\Om \varphi_{\ho} (Eu)\dx
\end{equation}
in $\wup\orn$; moreover, if $\varphi$ is convex and $p=1$ then
$$
F_{\ho}(u,\Om)=\int _\Om \varphi_{\ho} ({\mathcal E}u(x))\dx
+\int _\Om \varphi _{\ho} ^\infty \Bigl (\frac{dE^s u}{d\vert
  E^s u\vert }\Bigr )d\vert E^s u\vert
$$ 
in $\bd$, where $\varphi_{\ho}$ is described by an asymptotic formula.
However, in the second case, the question about the integral 
representation of the $\Gamma$-limit without the convexity assumption on $\varphi $ 
remained open. Notice that such an assumption is too strong in the 
vectorial calculus of variations. As shown in \cite{AE}, (see the proof of Theorem 5.1), 
the $\Gamma$-limit verifies the properties (1)-(3) and the invariance 
with respect to rigid motions follows from the fact that the 
approximating functionals $F_\e$ have this property. So we can apply 
Theorem \ref{integ} to obtain that
\begin{equation}\label{int2}
F_{\ho}(u,A)=\int_A f(x,\E u)dx +\int_{J_u\cap A} g(x, [u],\nu_u)d\hn\,,
\end{equation}
for every $u\in \sbd$ and every $A\in\omreg$. 
Now, from the integral representation (\ref{intw}) in $W^{1,1}(\Om,\R^N)$ and the 
relaxation theorem 3.5 in \cite{BFT}, one can easily see that
$$
f(x,\xi)=\varphi _{\ho}(\xi)\qquad\mbox{and}\qquad 
g(x, a,\nu)= \varphi _{\ho} ^\infty(a\odot\nu)
$$
for every $x\in\Om,\, a\in\R^N,\, \nu\in S^{N-1}$, and for every $\xi\in\msim.$
Notice that, since $\varphi _{\ho}$ is symmetric quasiconvex, that is 
$$ 
\varphi_{\ho}(\xi)\leq\intcauchy _A \varphi _{\ho}(\xi+{\mathcal E}\psi (x))dx
$$
for every $\psi \in W^{1,\infty }_0(A,\R ^N )$, $\xi\in\msim$ and for 
every bounded open subset $A$ of $\R ^N $, then  
$\varphi _{\ho} ^\infty$ is well defined.\\

{\bf Acknowledgements. } The authors wish to thank G. Dal Maso for many useful discussions concerning 
the subject of this paper, in particular for suggesting the example in Remark~\ref{rigid}. 
The work of Rodica Toader is part of the European Research Training Network ``Homogenization and Multiple Scales'' under contract
HPRN-2000-00109.


\begin{thebibliography}{99}

\bibitem{ACDM}{Ambrosio L., A. Coscia and G. Dal Maso}: Fine properties of 
functions with bounded deformation. {\it Arch. Rat. Mech. Anal.} {\bf 139} 
(1997), 201-238.  

\bibitem{ABC}{Ansini N., A. Braides and V. Chiad\`o Piat}: Homogenization of 
periodic multi-dimensional structures. {\it Boll. Un. Mat. Ital.} (8) {\bf 2-B} (1999), 735-758.

\bibitem{AE}{Ansini N. and F.\,B. Ebobisse}: Homogenization of periodic 
multi-dimensional structures II: the linearly elastic/perfect plastic case. 
{\it Adv. Math. Sci. Appl.} To appear.

\bibitem{AG}{Anzellotti G. and M. Giaquinta}: Existence of the displacement 
field for an elasto-plastic body subject to Henky's law and Von Mises' yield 
condition. {\it Manuscripta Math.} {\bf 32} (1980), 101-131.

\bibitem{BFT}{Barroso A.\,C., I. Fonseca and R. Toader}: A relaxation theorem 
in the space of functions of bounded deformation. {\it Ann. Sc. Norm. Sup. 
Pisa} {\bf 29} (2000), 19-49.

\bibitem{BCDM}{Bellettini G., A. Coscia and G. Dal Maso}: Special functions of 
bounded deformation. {\it Math. Z.} {\bf 228} (1998), 337-351.  

\bibitem{BFM}{Bouchitt\'e G., I. Fonseca and L. Mascarenhas}: A global method 
for relaxation. {\it Arch. Rat. Mech. Anal.} {\bf 145} (1998), 45-68.





\bibitem{DA}{De Giorgi E. and L. Ambrosio}: Un nuovo tipo di funzionale del 
calcolo delle variazioni. {\it Atti Accad. Naz. Lincei Rend. Cl. Sci. Mat. 
Natur.} {\bf 82} (1988), 199-210.


\bibitem{EBO2}{Ebobisse F.}: On lower semicontinuity of integral functionals 
in $LD(\Om )$. {\it Ric. di Matematica} {\bf XLIX} (2000), 65-76.




\bibitem{KH}{Kohn R.\,V.}: {\it New Estimates for Deformations in terms of 
their Strains}, Ph.D. Thesis, Princeton University, 1979.

\bibitem{MSC}{Matthies H., G. Strang and E. Christiansen}: The Saddle Point 
of a Differential Program. {\it Energy Methods in Finite Element Analysis},
Wiley, New York, 1979.

\bibitem{SUQ}{Suquet P.\,M.}: Un espace fonctionnel pour les equations de la 
plasticit\'e. {\it Ann. Fac. Sci. Toulouse} {\bf 1} (1979), 77-87.

\bibitem{TEM}{Temam R.}: {\it Probl\`emes Math\'ematiques en Plasticit\'e}. 
Gauthier-Villars, 1983.

\end{thebibliography}
\end{document}